\input amstex
\documentstyle{amsppt}
\magnification=\magstep1

\define\>{\rightarrow}
\define\<{\leftarrow}
\define\[{\lbrack}
\define\]{\rbrack}
\redefine\o{\circ}

\define\al{\alpha}
\define\be{\beta}

\define\de{\delta}

\define\si{\sigma}

\define\ph{\varphi}

\define\ps{\psi}

\define\Th{\Theta}
\define\La{\Lambda}

\define\Ph{\Phi}

\define\Om{\Omega}

\define\row#1#2#3{#1_#2,\ldots,#1_#3}
\define\x{\times}

\topmatter
\title A Cohomology for Vector Valued Differential Forms   \endtitle
\author Peter W. Michor \\
       Hubert Schicketanz \endauthor
\affil Institut f\"ur Mathematik der Universit\"at Wien, Austria \\
Institut des Math\'ematique, Universit\'e de Li\`ege, Belgium \endaffil
\address{Institut f\"ur Mathematik, Universit\"at Wien,
Strudlhofgasse 4, A-1090 Wien, Austria}\endaddress
\address{Universit\'e de Li\`ege, Institut de Math\'ematique, Avenue
des Tilleuls, 15; B-4000 Li\`ege, Belgium}\endaddress
\abstract{A rather simple natural outer derivation of the graded
Lie algebra of all vector valued differential forms with the
Fr\"olicher-Nijenhuis bracket turns out to be a differential and
gives rise to a cohomology of the manifold, which is functorial
under local diffeomorphisms. This cohomology is determined as
the direct product of the de Rham cohomology space and the graded
Lie algebra of "traceless" vector valued differential forms,
equipped with a new natural differential concomitant as graded
Lie bracket. We find two graded Lie algebra structures on the
space of differential forms. Some
consequences and related results are also discussed.}\endabstract
\endtopmatter
\document

\heading 1. Notation \endheading

\subheading{1.1 The Fr\"olicher-Nijenhuis bracket} Let $M$
be a smooth manifold of dimension $m$ throughout the paper. 

We
consider the space $\Om(M;TM) = \bigoplus_{k=0}^m\,\Om^k(M;TM)$
of all tangent bundle valued differential forms on $M$. Below
$K$ and $L$ will be elements of $\Om(M;TM)$ of degree $k$ and
$\ell$, respectively.
It is well known that $\Om(M;TM)$ is
a graded Lie algebra with the so called {\it
Fr\"olicher-Nijenhuis bracket} 
$$[\quad,\quad]: \Om^k(M;TM) \x \Om^\ell(M;TM) \> \Om^{k+\ell}(M;TM).$$
For its definition, properties, and notation we refer to
\cite{Mi, 1987}.

\bigpagebreak
\subheading{1.2}
In the investigation of the Lie algebra cohomology of the graded
Lie algebra
$(\Om(M;TM),\;[\quad,\quad])$ in \cite{Sch,1988} the following
exterior graded derivation of degree 1 appeared:
$$\de:\Om^k(M;TM) \> \Om^{k+1}(M;TM)$$
Before its definition we need another operator. Let the {\it
contraction} or {\it trace} $c:\Om^k(M;TM) \> \Om^{k-1}(M)$  be
given by $c(\ph \otimes X) = i_X\ph$, linearly extended. 
We also put 
$$\bar c\mid\Om^k(M;TM) := \tsize \frac{(-1)^{k-1}}{m-k+1}\;c$$ 
for reasons which become clear in lemma 2.1 below.

Then let $\de(K) := (-1)^{k-1}dc(K)\wedge \bold I$,  where
$\bold I = Id_{TM} \in \Om^1(M;TM)$ is the generator of the
center of the Fr\"olicher--Nijenhuis algebra.
This operator is a derivation with respect to the
Fr\"olicher-Nijenhuis-bracket, so we have $\de ([K,L]) =
[\de(K),L] + (-1)^k [K,\de(L)]$; furthermore $\de \o \de = 0$.
These properties are proved in \cite{Sch,1988} and are
straightforward to check, using lemma 2.1 below.  

\heading 2. The cohomology space $H(\Om(M;TM), \de)$ \endheading

\proclaim{2.1. Lemma} The mapping $j := (\quad\wedge \bold I):
\Om^k(M) \> \Om^{k+1}(M;TM)$ is a right inverse for $\bar c$
and the following diagram commutes for $k<m$:
$$\CD 
\Om^{k-1}(M)  @>j>> \Om^k(M;TM) @>\bar c>> \Om^{k-1}(M)\\
@VdVV               @V{\frac 1{m-k+1}\;\de}VV  @Vd VV \\   
\Om^k(M)      @>j>> \Om^{k+1}(M;TM) @>\bar c>> \Om^k(M)
\endCD\tag a  $$
\endproclaim
\demo{Proof} Both operators are local and in a coordinate system
we write $\bold I = Id_{TM} = \sum dx^i\otimes \partial_i$. Let $\ph \in
\Om^k(M)$, then we have 
$$\split c(\ph \wedge \bold I) &= c\left(\tsize\sum_i\ph
\wedge dx^i \otimes \partial_i\right) \\  
&= \tsize\sum_i\left(i_{\partial_i}\ph \wedge dx^i +
(-1)^k \ph \wedge i_{\partial_i}dx^i\right) \\  
&= (-1)^{k-1}k\ph + (-1)^km\ph = 
(-1)^k(m-k)\ph. \endsplit $$
The rest is a consequence of this.\qed\enddemo
\bigpagebreak
\subheading{2.2} Since $\bar c \o j = Id$,
the mapping $P:= j\o \bar c :\Om(M;TM) \>
\Om(M;TM)$ is a projection, so $P\o P = P$ and 
$$\split \Om(M;TM) &= \operatorname{im}P \oplus \operatorname{ker}P,\\
K &= \bar c(K)\wedge \bold I + K'.\endsplit \tag a $$
Note that $j$ is injective; this has the following consequences:
$K'$ in the decomposition (a) is characterized by $\bar
c(K') =0$. Moreover $\de(K) = 0$ if and only if $d\bar c(K) = 0$.
Finally $\operatorname{ker}P = \operatorname{ker}\bar c =
\operatorname{ker}c$; we will denote this space by 
$$\Cal C(M) = \bigoplus_{k=0}^m \Cal C^k(M). \tag b $$
Then $\Cal C^k(M) = \{\,K\in \Om^k(M;TM): \bar c(K) = 0 \,\} =
C^\infty(E^k(M))$ is the space of smooth sections of a certain natural
vector bundle over $M$, namely the subbundle
$\operatorname{ker}c$ of $\La^kT^*M \otimes TM$.
If $(U,x)$ is a chart on $M$,
the sections $dx^{i_1} \wedge \ldots \wedge dx^{i_k} \otimes
\partial_j$ with $i_1 < \cdots < i_k$ and $j \ne i_l$ for all $l$
give a local framing of this bundle. Note that $\Cal C^0(M) =
\Cal X(M)$, the space of all vector fields, and that $\Cal C^1(M)$ is
the space of all traceless endomorphisms $TM \> TM$. For this
reason elements of $\Cal C(M)$ will be called {\it traceless}
vector valued differential forms. Note that $\Cal C^m(M) = 0$
since $\bar c:\Om^m(M;TM) \> \Om^{m-1}(M)$ is a linear isomorphism.  

\subheading{2.3} Let us define the natural bilinear concomitant
$$S:\Om^k(M;TM) \x \Om^\ell(M;TM) \> \Om^{k+\ell-2}(M)$$
by
$S(\ph \otimes X,\ps\otimes Y) =  i_Y\ph \wedge i_X\ps$
for decomposable vector valued forms.

\proclaim{Lemma} Then we have
$$c([K,L]) = (-1)^k\Th(K)c(L) - (-1)^{k\ell}\Th(L)c(K) 
- (-1)^k dS(K,L).$$ 
\endproclaim

\demo{Proof} Since both sides are local in $K$ and $L$
we may assume that $K = \ph \otimes X$ and $L = \ps \otimes Y$
are decomposable. Then by formula \cite{Mi, 1987, 1.7.7} we have
$$\multline [K,L] = \ph \wedge \ps \otimes [X,Y] + \ph \wedge
\Th(X)\ps \otimes Y - \Th(Y)\ph \wedge \ps \otimes X + \\
+ (-1)^k\left(d\ph \wedge i_X\ps \otimes Y + i_Y\ph \wedge d\ps
\otimes X\right).\endmultline \tag a $$
This implies the lemma by a straightforward computation.
\qed\enddemo
\bigpagebreak

\proclaim{2.4. Theorem} $(\Cal C(M),\;[\quad,\quad]^c)$ is a
graded Lie algebra, where 
$$[K,L]^c := [K,L] - \tsize\frac{(-1)^\ell}{m-k-\ell+1}\; dS(K,L) \wedge
\bold I. \tag a $$ 
It is a quotient of a subalgebra of $(\Om(M;TM),\; [\quad,\quad])$.
The bracket $[\quad,\quad]^c$ is a natural bilinear differential
concomitant of order $1$, that means $f^*[K,L]^c =
[f^*K,f^*L]^c$ for each local diffeomorphism $f$.
\endproclaim
This theorem will be proved jointly with theorem 2.5 below.

\bigpagebreak
\proclaim{2.5. Theorem} The cohomology of the graded Lie
algebra  $\Om(M;TM)$ is decomposed into:
$$\split H^k(\Om(M;TM),\de) &\cong H^{k-1}_{dR}(M) \oplus
(\operatorname{ker}\bar c)^k \\
&= H^{k-1}_{dR}(M) \oplus \Cal C^k(M) \endsplit \tag a $$
The induced bracket 
$$[\quad,\quad]:H^k(\Om(M;TM),\de) \x H^\ell(\Om(M;TM),\de)  \>
H^{k+\ell}(\Om(M;TM),\de)$$
corresponds to the direct product of 
the graded Lie algebra $(\Cal C(M),\;[\quad,\quad]^c)$ with the
abelian algebra $(H^{*-1}_{dR}(M),0)$.   
\endproclaim

\demo{Proof of 2.4 and 2.5} By diagram 2.1.a we have induced
mappings in cohomology
$j^\sharp = (\quad\wedge \bold I)^\sharp: H^{*-1}_{dR}(M) \> 
H^*(\Om(M;TM),\de)$ 
and $\bar c^\sharp: H^*(\Om(M;TM),\de) \>
H^{*-1}_{dR}(M)$ ; again $\bar c^\sharp \o j^\sharp = Id$, so
$P^\sharp = j^\sharp \o \bar c^\sharp$ 
is also a projection in cohomology, and we have the
decomposition 
$$ H^k(\Om(M;TM),\de) \cong (\operatorname{im}P^\sharp)^k \oplus
(\operatorname{ker}P^\sharp)^k = H^{k-1}_{dR}(M) \oplus
\Cal C^k(M).$$
Here we used that $j^\sharp:H^{k-1}_{dR}(M) \> (\operatorname{im}P^\sharp)^k$
is a linear isomorphism, and that $\operatorname{ker}P^\sharp =
\operatorname{ker}P$ since $\de = \pm dc(\quad)\wedge \bold I$ implies
$\operatorname{im}\de \subset \operatorname{im}j =
\operatorname{im}P$ and $\operatorname{ker}P =
\operatorname{ker}\bar c \subset \operatorname{ker}\de$.

Since the differential $\de$ is a graded derivation,
$\operatorname{ker}\de$ is a graded Lie subalgebra of
$(\Om(M;TM),[\quad,\quad])$, the space
$\operatorname{im}\de$ is an ideal in $\operatorname{ker}\de$
and the cohomology space $H(\Om(M;TM),\de)$ is a graded Lie
algebra. 

It remains to investigate the induced bracket. 
Let $K = \bar c(K) \wedge \bold I  + K' \in
(\operatorname{ker}\de)^k \subset 
\Om^k(M;TM) = (\operatorname{im}P)^k \oplus
(\operatorname{ker}P)^k$ and similarly 
$L = \bar c(L) \wedge \bold I + L'$, then by using formula 
\cite{Mi, 1987, 1.7.6} we may compute as follows: 
$$\allowdisplaybreaks\split
[K,L] &= [\bar c(K) \wedge \bold I,\bar c(L) \wedge \bold I] +
[\bar c(K) \wedge \bold I,L'] + [K',\bar c(L) \wedge \bold I] +
[K',L']\\
&= 0 -(-1)^{k\ell}\Th(\bar c(L) \wedge \bold I)\bar c(K) \wedge
\bold I + (-1)^kd\bar c(K) \wedge i_{\bold I}(\bar c(L) \wedge
\bold I) \\
&\quad + 0 -(-1)^{k\ell}\Th(L')\bar c(K) \wedge \bold I +
(-1)^k d\bar c(K) \wedge i_{\bold I}L' \\
&\quad + 0 + \Th(K')\bar c(L) \wedge \bold I -
(-1)^{k\ell+\ell}d\bar c(L) \wedge i_{\bold I}K' \\
&\quad +[K',L'] \\
&= \Th(K')\bar c(L) \wedge  \bold I - (-1)^{k\ell}\Th(L')\bar
c(K) \wedge \bold I + [K',L'], 
\endsplit $$ 
since by the definition of $\de$ we have $d\bar c(K) =0$ and
$d\bar c(L) = 0$ for $k$,~$\ell < m$. By 2.2 we have $\bar c(K')
= 0$ and $\bar c(L') = 0$, so by 2.2, 2.3, and 2.4 we have
$$\split [K',L'] &= \bar c([K',L']) \wedge \bold I + [K',L']^c \\
&= \tsize \frac{(-1)^\ell}{m-k-\ell+1} dS(K',L') \wedge \bold I +
[K',L']^c.\endsplit \tag b $$ 
Since $\bar c(K)$ is closed, $\Th(L')\bar c(K) = i(L')d\bar c(K)
- (_1)^{\ell -1}di(L')\bar c(K)$ is exact, likewise 
$\Th(K')\bar c(L)$ is exact. Therefore in 
$H(\Om(M;TM),\de) = H^{*-1}(M) \oplus \Cal C(M)$ we have
$[\al + K',\be + L'] = [K',L']^c$.
So $H^{*-1}(M)$ is an abelian ideal of  $H^*(\Om(M;TM),\de)$,
the subspace $(\Cal C(M),[\quad,\quad]^c)$ is an ideal of
$H^*(\Om(M;TM),\de)$, and a quotient of the graded Lie algebra 
$j(Z^{*-1}(M)) \oplus \Cal C(M)$. The differential concomitant
$[\quad,\quad]^c$ is natural since its composants (in 2.4) are
all natural.
\qed\enddemo 
\bigpagebreak

\heading 3. Extensions of the graded Lie algebra $\Cal C(M)$ \endheading

\subheading{3.1. Graded Lie subalgebras of $\Om(M;TM)$} Let 
$$\alignat 2
\Om^{*-1}(M) &= \tsize\bigoplus_{k=0}^{m-1}\Om^k(M), &\qquad
Z^{*-1}(M) &= \tsize\bigoplus_{k=0}^{m-1}Z^k(M),\\
B^{*-1}(M) &= \tsize\bigoplus_{k=0}^{m-1}B^k(M), &\qquad
H^{*-1}_{dR}(M) &= \tsize\bigoplus_{k=0}^{m-1}H^k_{dR}(M)
\endalignat $$
be the graded spaces of
de~Rham forms, cycles, boundaries, and cohomology classes,
respectively, where we exchanged 
the top degree spaces for $0$. 

\proclaim{Theorem} With this notation we have:
\roster
\item $B^{*-1}(M) \wedge \bold I$ is an abelian ideal in $\Om^*(M;TM)$.
\item $Z^{*-1}(M) \wedge \bold I$ is an abelian ideal in $\Om^*(M;TM)$.
\item $\Om^{*-1}(M) \wedge \bold I$ is a graded Lie subalgebra
of $\Om^*(M;TM)$. 

Therefore $(\Om^{*-1}(M),\; [\quad,\quad])$ is a
graded Lie algebra, where the bracket induced from the
Fr\"olicher-Nijenhuis one looks as follows for $\ph \in
\Om^k(M)$ and $\ps \in \Om^\ell(M)$:
$$\split[\ph,\ps] &= (-1)^{k-1}\left(d\ph \wedge \ell \ps
-(-1)^{(k-1)(\ell-1)}d\ps \wedge k\ph\right) \\
&= (-1)^{k-1}\left(\Th(\bold I)\ph \wedge i_{\bold I} \ps
-(-1)^{(k-1)(\ell-1)}\Th(\bold I)\ps \wedge i_{\bold I}\ph\right) 
\endsplit $$
\endroster
\endproclaim
\demo{Proof} Straightforward computations using \cite{Mi, 1987, 1.7}.
\qed\enddemo

\subheading{3.2. Extensions of $\Cal C(M)$} If we collect all
relevant parts in the proof  of theorem 2.5 we get  for  
$K = \bar c(K) \wedge  \bold I + K'\in(Z^{k-1}(M) \wedge \bold
I) \oplus \Cal C^k(M)$ and 
$L = \bar c(L) \wedge \bold I + L' \in (Z^{\ell-1}(M) \wedge
\bold I) \oplus \Cal C^\ell(M)$ the following
formula: 
$$ [K,L] =  \left(\Th(K')\bar c(L) -(-1)^{k\ell}\Th(L')\bar c(K)
+ \si(K'L')\right) \wedge \bold I + [K'L']^c, \tag a $$
where $\si:\Cal C^k(M) \x \Cal C^\ell(M) \> B^{k+\ell-1}(M)
\subset Z^{k+\ell-1}(M)$ is defined by 
$$\si(K',L') := \tsize\frac{(-1)^\ell}{m-k-\ell+1} dS(K',L').$$ 
Now by 3.1.2 we have an exact sequence of graded Lie algebras 
$$ 0 \> Z^{*-1}(M) @>j>> (Z^{*-1}(M) \wedge \bold I) \oplus \Cal C(M)
\> \Cal C(M) \> 0, \tag b $$
which describes an abelian extension of $\Cal C(M)$. By 2.4 we
have 
$$\split\Th([K',L']^c) &= \Th([K',L'] - \si(K',L')\wedge \bold I)\\
&= [\Th(K'),\Th(L')] -\si(K',L') \wedge \Th(\bold I)|Z(M) \\
&= [\Th(K'),\Th(L')]. \endsplit \tag c $$ 
So $\Th$ gives a graded Lie module structure over $\Cal C(M)$ to
$Z^{*-1}(M)$ and consequently $\si$ is
a cocycle for the graded Lie algebra cohomology of $(\Cal
C(M),\;[\quad,\quad]^c)$ with coefficients in the $\Cal
C(M)$-module $Z^{*-1}(M)$, because
the second cohomology classifies equivalence classes of abelian
extensions, just as in the non graded case. For the convenience
of the reader we sketch this in 3.3 below. A direct check shows
that indeed the cocycle equation for $\si$ is valid and is
equivalent to the graded Jacobi identity for $(\Cal
C(M),\;[\quad,\quad]^c)$. 

\bigpagebreak
\subheading{3.3. Cohomology of graded Lie algebras} In the
following we write down the definitions for the graded
cohomology of a $\bold Z$-graded Lie algebra $L = \bigoplus_{k
\in \bold Z}L^k$ with coefficients in a graded $L$-module $V$.
So $V = \bigoplus_{k \in \bold Z}V^k$ is a graded vector space and
$\Th: L \> \operatorname{End}(V)$ is a homomorphism of graded
algebras (of degree 0, where the bracket on
$\operatorname{End}(V)$ is the graded commutator). Our
definitions are different from but equivalent to those of
\cite{Le}, which are also used in \cite{Sch}; we obey Quillen's
rule strictly.

Let $\La^p(L;V)^q$ be the space of all $p$-linear mappings
$\Ph:L \x \cdots \x L \> V$, which are of degree $q$, i.e.,
$\Ph(\row X1p) \in V^{x_1+\cdots+x_p+q}$ for $X_i \in L^{x_i}$,
and which are alternating in the sense that $\Ph$ reacts to
interchanging $X_i$ and $X_{i+1}$ with the sign $-(-1)^{x_ix_{i+1}}$.

Now let us define the differential 
$\partial:\La^p(L;V)^q \> \La^{p+1}(L;V)^q$ by 
$$\multline (\partial\Ph)(\row X0p) 
= \sum_i (-1)^{\al_i+x_iq}\,\Th(X_i)\Ph(X_0,\ldots,
\overbrace{X_i},\ldots,X_p) \\
 + \sum_{i<j}(-1)^{\al_i+\al_j-x_ix_j}\Ph([X_i,X_j],X_0,\ldots,
\overbrace{X_i},\ldots,\overbrace{X_j},\ldots X_p),
\endmultline $$
where $\al_i = i+x_i(x_0+\dots+x_{i-1})$, and where the brace
over a symbol means that it has to be deleted. Then one may
check that $\partial \o \partial = 0$, and we denote by 
$H^p(L;V)^q := H^p(\La^*(L;V)^q,\partial)$ the resulting
cohomology and call it the {\it graded cohomology} of the graded
Lie algebra $L$.

\proclaim{Theorem} $H^2(L;V)^0$ is isomorphic to the set of
equivalence classes of abelian extensions of $L$ by $V$.
\endproclaim
The proof of this is completely analogous to the non-graded case
after the insertion of some obvious signs.

\bigpagebreak
\subheading{3.4. The Nijenhuis-Richardson bracket} Recall from
\cite{Ni-Ri, 1967} or \cite{Mi, 1987}, that there is a natural
graded Lie algebra structure on $\Om^{*-1}(M;TM)$, given by 
$$[K,L]^{\wedge} = i(K)L -(-1)^{(k-1)(\ell-1)}i(L)K.$$

\proclaim{Theorem} For the Nijenhuis-Richardson bracket we have 
\roster
\item $\Cal C^{*-1}(M)$ is a graded Lie subalgebra of
$(\Om^{*-1}(M;TM), [\quad,\quad]^\wedge)$. 
\item $\Om^*(M) \wedge \bold I$ is a graded Lie
subalgebra of $(\Om^{*-1}(M;TM), [\quad,\quad]^\wedge)$. Therefore
$(\Om^*(M),\;[\quad,\quad]^\wedge)$ is a graded Lie algebra,
where the induced bracket looks as follows for $\ph \in
\Om^k(M)$ and $\ps \in \Om^\ell(M)$:
$$\split [\ph,\ps]^\wedge &= (\ell-k)\ph\wedge\ps \\
&= \ph \wedge i_{\bold I}\ps - (-1)^{k\ell}\ps \wedge i_{\bold I}\ph\qed
\endsplit $$
\endroster
\endproclaim

None of these two subalgebras is an ideal, so there is no
extension. For the structure of the whole algebra in terms of
the subalgebras see \cite{Mi, 1988}

Compare the formulas here and in 3.1. They rise the
question, whether there is general procedure behind them.

\enddocument